\documentclass[10pt]{amsart}
\usepackage{amssymb, amsthm, amsmath}
\usepackage{latexsym}
\usepackage{epsfig}
\usepackage{xy}
\usepackage{psfrag}

\theoremstyle{plain}
\newtheorem{prop}{Proposition}
\newtheorem{fact}{Fact}
\newtheorem{thm}{Theorem}

\newtheorem{lemma}{Lemma}
\theoremstyle{definition}
\newtheorem*{remark}{Remark}
\newtheorem*{example}{Example}

\def\Z{{\mathbb Z}}

\def\C{{\mathbb C}}

\def\l{{\lambda}}

\def\Sh{{\mathcal S}}

\def\s{{\sigma}}
\def\t{{\tau}}

\newcommand{\uu}{{\bf u}}
\newcommand{\vv}{{\bf v}}
\newcommand{\w}{{\bf w}}
\newcommand{\x}{{\bf x}}
\newcommand{\y}{{\bf y}}
\newcommand{\z}{{\bf z}}

\newcommand{\gequ}{\geqslant}
\newcommand{\lequ}{\leqslant}

\newcommand{\noin}{\noindent}
\newcommand{\wt}{\widetilde}

\newcommand{\wh}{\widehat}

\newcommand{\hole}{{\mathrm o}}

\newcommand{\reflemma}[1]{Lemma~\ref{#1}}

\newcommand{\picB}[1]{\includegraphics[scale=0.75]{#1}}
\newcommand{\picC}[1]{\includegraphics[scale=0.60]{#1}}

\newcommand{\p}{\negmedspace + \negmedspace}
\newcommand{\m}{\negmedspace - \negmedspace}

\newenvironment{abcenum}{\begin{enumerate}}{\end{enumerate}}

\psfrag{a}{$a$}
\psfrag{b}{$b$}
\psfrag{x}{$x$}
\psfrag{y}{$y$}
\psfrag{i}{$i$}
\psfrag{i'}{$i'$}
\psfrag{uu}{$\bf u$}
\psfrag{vv}{$\bf v$}
\psfrag{xx}{$\bf x$}
\psfrag{yy}{$\bf y$}
\psfrag{zz}{$\bf z$}
\psfrag{1}{$1$}
\psfrag{2}{$2$}
\psfrag{3}{$3$}
\psfrag{4}{$4$}
\psfrag{1'}{$1'$}
\psfrag{2'}{$2'$}
\psfrag{3'}{$3'$}
\psfrag{4'}{$4'$}
\psfrag{H}{$\hole$}
\psfrag{H'}{$\hole'$}
\psfrag{H2}{$\hole_2$}

\begin{document}

\title[Littlewood-Richardson rules]
{Littlewood-Richardson rules for Grassmannians}
\author{Anders Skovsted Buch, Andrew Kresch, and Harry Tamvakis}
\date{June 22, 2003}
\thanks{The authors were supported in part by NSF Grant DMS-0070479
  (Buch), an NSF Postdoctoral Research Fellowship (Kresch),
  and NSF Grant DMS-0296023 (Tamvakis).}
\address{Matematisk Institut, Aarhus Universitet, Ny Munkegade, 8000
  {\AA}rhus C, Denmark}
\email{abuch@imf.au.dk}
\address{Department of Mathematics, University of Pennsylvania,
209 South 33rd Street,
Philadelphia, PA 19104-6395, USA}
\email{kresch@math.upenn.edu}
\address{Department of Mathematics, Brandeis University - MS 050,
P. O. Box 9110, Waltham, MA
02454-9110, USA}
\email{harryt@brandeis.edu}

\maketitle

\section{Introduction}
\label{intro}

The classical Littlewood-Richardson rule \cite{LR} describes the
structure constants obtained when the cup product of two Schubert
classes in the cohomology ring of a complex Grassmannian is written as
a linear combination of Schubert classes.  It also gives a rule for
decomposing the tensor product of two irreducible polynomial
representations of the general linear group into irreducibles, or
equivalently, for expanding the product of two Schur $S$-functions in
the basis of Schur $S$-functions.  In this paper we give a short and
self-contained argument which shows that this rule is a direct
consequence of Pieri's formula \cite{P} for the product of a Schubert
class with a special Schubert class.

There is an analogous Littlewood-Richardson rule for 
the Grassmannians which parametrize maximal isotropic subspaces
of $\C^n$, equipped with a symplectic or orthogonal form. The precise
formulation of this rule is due to Stembridge \cite{St}, working in
the context of Schur's $Q$-functions \cite{S2}; the connection to 
geometry was shown by Hiller and Boe \cite{HB} and Pragacz \cite{Pr}.
The argument here for the type $A$ rule works equally well in these more
difficult cases and gives a simple derivation of Stembridge's rule
from the Pieri formula of \cite{HB}.

Currently there are many proofs available for the classical
Littlewood-Richardson rule, some of them quite short.  The proof of
Remmel and Shimozono \cite{RS} is also based on the Pieri rule; see
the recent survey of van Leeuwen \cite{vL} for alternatives.  In
contrast, we know of only two prior approaches to Stembridge's rule
(described in \cite{St,HH} and \cite{Sh}, respectively), both of which
are rather involved.

The argument presented here proceeds by defining an abelian group
${\mathbb H}$ with a basis of Schubert symbols, and a bilinear product
on ${\mathbb H}$ with structure constants coming from the
Littlewood-Richardson rule in each case. Since this rule is compatible
with the Pieri products, it suffices to show that ${\mathbb H}$ is an
{\em associative} algebra.  The proof of associativity is based on
Sch\"utzenberger slides in type $A$, and uses the more general slides
for marked shifted tableaux due to Worley \cite{W} and Sagan \cite{Sa}
in the other Lie types.  In each case, we need only basic properties of
these operations which are easily verified from the definitions.  Our
paper is self-contained, once the Pieri rules are granted.

The work on this article was completed during a fruitful visit to
the Mathematisches Forschungsinstitut Oberwolfach, as part
of the Research in Pairs program. It is a pleasure to thank the 
Institut for its hospitality and stimulating atmosphere. 
We also thank the referee for a careful reading of our paper and 
for some well-placed final touches to the exposition.

\section{The Littlewood-Richardson rule for type $A$ Grassmannians}
\label{lrsec}

Let $X=G(k,n)$ be the Grassmannian of $k$-dimensional linear subspaces
of $\C^n$ and set $m=n-k$. For each partition $\l$ whose Young diagram
is contained in the $k\times m$ rectangle $(m^k)$, there is a {\em
Schubert class} $\sigma_\lambda$ in the cohomology ring $H^*(X,\Z)$.
If a partition $\wt \lambda \subset (m^k)$ can be obtained from
$\lambda$ by adding a horizontal strip with $p$ boxes, then we write
$\lambda \xrightarrow{p} \wt{\lambda}$.  The Pieri rule \cite{P}
states that for each $p\lequ m$, $\sigma_p \cdot \sigma_\lambda$ is
equal to the sum of all $\s_{\wt{\lambda}}$ for which $\lambda
\xrightarrow{p} \wt\lambda$.

In this section, we will prove that the Littlewood-Richardson rule
holds in the ring $H^*X=H^*(X,\Z)$. We note however that the argument
requires only two facts about this ring: (i) the classes $\sigma_\l$
for $\l\subset (m^k)$ form a basis of $H^*X$, and (ii) the Pieri rule
holds in $H^*X$. An easy induction shows that the special Schubert
classes $\sigma_p$ for $1\lequ p\lequ m$ generate the entire ring
$H^*X$. This also follows from the Giambelli formula, which is a
direct consequence of Pieri's rule.  Let $\lambda^\vee = (m-\lambda_k,
\dots, m-\lambda_1)$ denote the dual partition of $\lambda$.

A {\em tableau} $T$ of skew shape $\l/\mu$ is a filling of the boxes of
$\l/\mu$ with positive integers such that the entries are weakly 
increasing along each row and strictly increasing down each column. 
The {\em content} of $T$ is the sequence whose $i$th element is 
the number of boxes of $T$ containing $i$. 
The {\em word} $\w=\w(T)$ of $T$ is the sequence obtained by reading the
entries of $T$ going from right to left in successive rows, starting
with the top row. We say
that $\w=w_1\ldots w_r$ is a {\em lattice word} and 
that $T$ is a 
{\em Littlewood-Richardson tableau} (or {\em LR tableau}) if 
the number of occurrences of $i$ among  $w_1\ldots w_j$ 
is not less than the number of occurrences of $i+1$, 
for all $i$ and $j$ with $1\lequ j\lequ r$.

Given three partitions $\lambda, \mu, \nu \subset (m^k)$, define
$c(\lambda, \mu ;\nu)$ to be the number of LR 
tableaux of shape $\l^\vee/\mu$ with content $\nu^\vee$.  (If
$\mu$ is not contained in $\l^\vee$ then we set 
$c(\l, \mu;\nu) = 0$.)

\begin{prop}
\label{assprop}   For any three partitions 
  $\l,\mu,\nu \subset (m^k)$ and integer
  $p \lequ m$, we have
\begin{equation}
\label{asseq}
\sum_{\l\xrightarrow{p} 
\widetilde\l} c(\widetilde\lambda, \mu;\nu)
   = \sum_{\mu \xrightarrow{p} 
\widetilde \mu} c(\lambda, \widetilde\mu;\nu) \,.
\end{equation}
\end{prop}
\begin{proof}
  We assume here familiarity with Sch\"utzenberger's jeu de
  taquin (explained e.g.\ in \cite[\S 1.2]{F}).
  Given a skew tableau $T$ and an empty box   which is an 
  inner corner of $T$, we may perform Sch{\"u}tzenberger slides
  to obtain a new skew tableaux $T'$; the empty box slides to an outer 
corner of $T$.

\begin{fact}
$T$ is an LR tableau if and only if $T'$ is an LR tableau.
\end{fact}

This follows immediately from the definitions; alternatively, 
it is a consequence of the well-known fact that
plactic relations on words preserve the lattice property.
For the direct
implication, it suffices to consider a single vertical slide as
displayed below. In the figure, the symbols $\uu$, $\x$, $\y$, 
$\z$ and $\vv$ 
denote the words of their respective subsets in the tableau. In 
particular, they are read from right to left.
\[ \picB{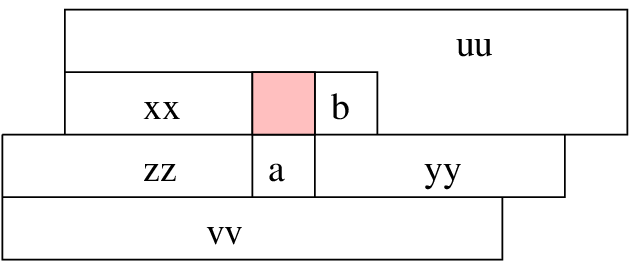} 
   \raisebox{1cm}{\hspace{10pt}$\longmapsto$\hspace{10pt}}
   \picB{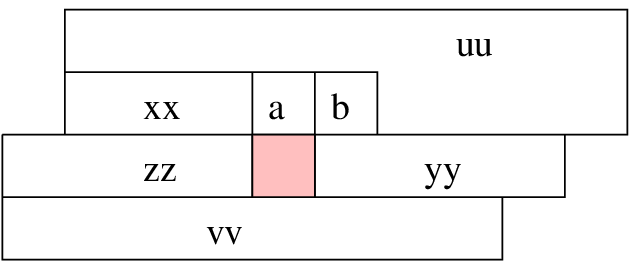} 
\]
We must check that the word $\uu b a\x\y\z\vv$ of the resulting
tableau is a lattice word.  This is true because $\uu b \x \y a \z
\vv$ is a lattice word, and the tableau inequalities imply that there
are at least as many $a$'s in the word $\z$ as there are $(a-1)$'s in
$\x$.  A similar argument shows that reverse slides also preserve the
lattice property.

We shall call an empty box contained inside the skew shape
$\l^\vee/\mu$ a {\em hole}.  Given an LR tableau on a shape
$\l^\vee/\widetilde\mu$ such that $\mu \xrightarrow{p} \widetilde\mu$,
we can use Sch{\"u}tzenberger slides starting from the holes contained
in $\widetilde\mu/\mu$, in right to left order, to obtain another LR
tableau of some shape $\widetilde\l^\vee/\mu$.  Define the {\em
  sliding path} of each such hole to be the set of boxes it occupies
during the sliding process.
\[ \picC{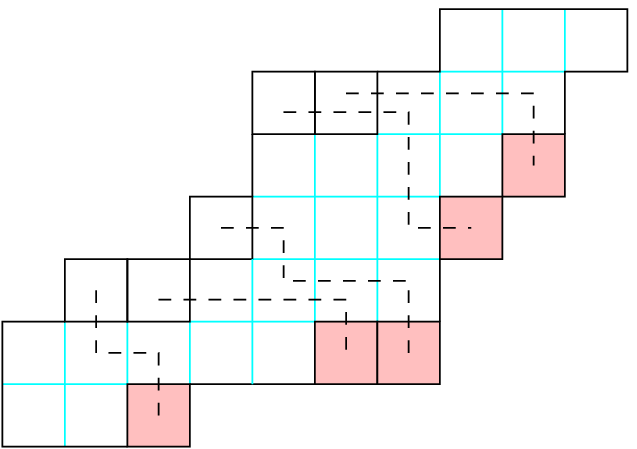} \]

\begin{fact}
Two distinct sliding paths cannot cross each other.
\end{fact}

More precisely, if a hole is at a given position during its slide, then the 
boxes in any subsequent sliding path must all lie strictly left or weakly
below that position. For otherwise, at some point a hole will slide right 
to occupy the position vacated by a vertical slide in the previous sliding
path. Depicting the vertical slide as
\[ \picB{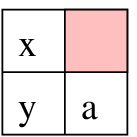}
   \hspace{10pt} \raisebox{.42cm}{$\longmapsto$} \hspace{10pt}
   \picB{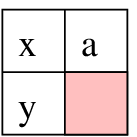}
\]
we must have $y \lequ a$, and hence a subsequent hole, having arrived at
position $x$, will slide down to position $y$.
Since different sliding paths
cannot cross each other, it follows that $\l \xrightarrow{p} \wt{\l}$.
Furthermore the entire process can be inverted using reverse slides.
This gives a bijective proof of identity (\ref{asseq}).
\end{proof}

\medskip

The following theorem is one out of many equivalent statements of the
classical Little\-wood-Richard\-son rule.
 
\begin{thm}  
\label{typeathm}
The constant $c(\lambda, \mu; \nu)$ is the coefficient of
  $\sigma_\nu$ in the product $\sigma_\lambda \cdot \sigma_\mu$.
\end{thm}
\begin{proof}
  Let ${\mathbb H}$ be the free abelian group generated by symbols
  $s_\lambda$ for all partitions $\lambda \subset (m^k)$.  We define a
  bilinear operator ``$\circ$'' on ${\mathbb H}$ by
\[ s_\lambda \circ s_\mu = \sum_\nu c(\lambda,\mu;\nu) \, s_\nu \,. \]
The operator $\circ$ is, a priori, neither commutative nor
associative.

It is easy to see that there is a unique LR tableau of shape $\l^\vee$
and a unique LR tableau of shape $(m^k)/\mu$, and that these tableaux
have contents $\l^\vee$ and $\mu^\vee$, respectively. It follows
that $s_{\emptyset}$ acts as a
left and right identity in ${\mathbb H}$.  By taking $\lambda =
\emptyset$ in Proposition \ref{assprop} we deduce that $s_p\circ s_\mu
= \sum s_{\wt\mu}$ where the sum is over $\mu \xrightarrow{p} \wt\mu$.
Similarly one obtains $s_\lambda \circ s_p = \sum s_{\wt\lambda}$ by
setting $\mu = \emptyset$; in other words, the operator $\circ$
satisfies the Pieri rule.

Equation (\ref{asseq}) is therefore equivalent to the associativity
relation $(s_\lambda \circ s_p) \circ s_\mu = s_\lambda \circ (s_p
\circ s_\mu)$. It follows that the elements $s_p$ for $1\lequ p\lequ
m$ generate an associative subalgebra of ${\mathbb H}$.  Using the
same Pieri induction as before, one sees that this subalgebra is the
entire algebra ${\mathbb H}$. We conclude that the linear map $H^*X
\to {\mathbb H}$ given by $\sigma_\lambda \mapsto s_\lambda$ is an
isomorphism of (associative) rings.
\end{proof}

\begin{remark} {\bf 1)}
  In its usual formulation, the Littlewood-Richardson rule states that
  the coefficient $c(\l,\mu;\nu)$ is equal to the number of LR
  tableaux of shape $\nu/\l$ with content $\mu$. To see this, note
  that the identity $c(\l,\mu; (m^k))=\delta_{\l,\mu^\vee}$ holds by
  definition (this corresponds to Poincar\'e duality in $H^*X$).
  It follows that
\[
c(\l,\mu;\nu)\s_{(m^k)} = \s_{\nu^\vee} (\s_\l \s_\mu) =
(\s_{\nu^\vee} \s_\l) \s_\mu = c(\nu^\vee,\l; \mu^\vee) \s_{(m^k)}
\]
and hence $c(\l,\mu;\nu)=c(\nu^\vee,\l; \mu^\vee)$, as required.
Alternatively, a bijective proof of this equality may be obtained 
using \cite[Prop. 5.1.2]{F}. 

\medskip
\noin {\bf 2)} The above argument may be applied to derive other forms
of the Littlewood-Richardson rule. For example, it gives a short
proof of the puzzle rule of Knutson, Tao and Woodward
\cite{KTW}. In the language of puzzles, Sch\"utzenberger slides
correspond to a subset of the propagations described in \cite{KT}
(those which involve only non-equivariant puzzle pieces).

\end{remark}

\section{The Littlewood-Richardson-Stembridge 
rule for maximal isotropic Grassmannians}

The odd orthogonal Grassmannian $Y=OG(n,2n+1)$ parametrizes
$n$-dimen\-sional isotropic linear subspaces of $\C^{2n+1}$ with respect
to a nondegenerate orthogonal form. The cohomology ring $H^*(Y,\Z)$
has a basis of Schubert classes $\t_{\l}$, indexed by strict
partitions $\l$ (i.e.\ with distinct parts) such that
$\l\subset\rho_n$, where $\rho_n=(n,n-1,\ldots,1)$.  For each strict
$\l\subset\rho_n$, define $\l^{\vee}\subset \rho_n$ as the strict
partition whose parts complement the parts of $\l$ in the set
$\{1,\ldots,n\}$. The {\em shifted diagram} $\Sh(\l)$ is obtained from
the Young diagram of $\l$ by indenting the $i$th row by $i-1$ columns,
for each $i\gequ 1$. For skew diagrams we set
$\Sh(\l/\mu)=\Sh(\l)\smallsetminus\Sh(\mu)$.  For example, if $n = 7$,
$\lambda = (5,3,1)$, and $\mu = (5,2)$ then $\Sh(\l^\vee/\mu)$ is the
diagram:
\[ \picC{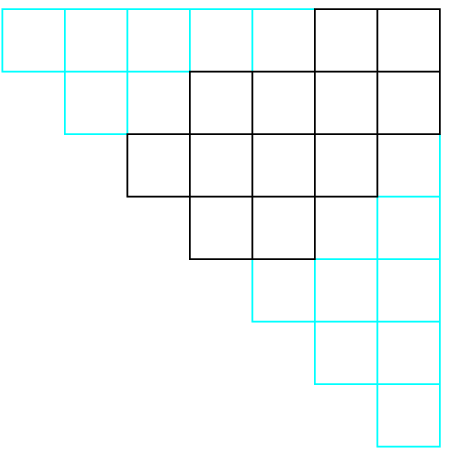} \]

Recall that a {\em border strip} is an edge-connected skew diagram
that contains no $2\times 2$ block of squares.  As before, we write
$\lambda \xrightarrow{p} \wt{\lambda}$ if the partition $\wt \lambda
\subset \rho_n$ can be obtained from $\lambda$ by adding a horizontal
strip of length $p$.  In this case, the shifted skew diagram
$\Sh(\wt{\l}/\l)$ is a union of border strips.  The Pieri rule for
$OG(n,2n+1)$, due to Hiller and Boe \cite{HB}, states that
\begin{equation}
\label{bpieri}
\t_p \cdot \t_\lambda=\sum 2^{N(\wt{\l}/\l)}\,\t_{\wt{\lambda}},
\end{equation}
where the sum is over strict 
$\widetilde\lambda \subset \rho_n$ with $\lambda \xrightarrow{p}
\widetilde\lambda$, and $N(\wt{\l}/\l)$ is one less than the number
of border strip components of $\Sh(\wt{\l}/\l)$. The Pieri rule implies
that the special Schubert classes $\t_p$ for $1\lequ p\lequ n$
generate $H^*(Y,\Z)$. 

Let $A$ be the ordered alphabet $1'<1<2'<2<\cdots$; the symbols
$1',2',\ldots$ are said to be {\em marked}.  A {\em shifted tableau}
$T$ on the shifted skew shape $\Sh(\l/\mu)$ is a filling of the boxes
of $\Sh(\l/\mu)$ with symbols from $A$ such that (i) the entries are
weakly increasing along each row and down each column, and (ii) each
row contains at most one $i'$ and each column contains at most one
$i$, for every integer $i\gequ 1$. The content of $T$ is the partition
whose $i$th part is the number of boxes with entry $i$ or $i'$ in $T$,
while the word $\w=\w(T)$ of $T$ is defined as in Section \ref{lrsec}.

For any integer $i$ we set $\wh{i'} = i$ and $\wh{i} = (i\p 1)'$.  If
$\w = w_1 w_2 \dots w_p$ is a word of marked and unmarked integers
$w_j$, then we write $\wh{\w} = \wh{w}_p \dots \wh{w}_2 \wh{w}_1$.  We
say that $\w$ is an {\em LRS word} if (i) $\w\wh{\w}$ is a {\em
lattice word}, i.e.\ every $i$ or $i'$ in $\w\wh{\w}$ is preceded by
more occurrences of $i-1$ than of $i$, for all $i$, and (ii) the last
occurrence of $i'$ in $\w$ (if any) is followed by at least one $i$,
for all $i\gequ 1$. A tableau $T$ is a {\em
Littlewood-Richardson-Stembridge tableau} (or {\em LRS tableau}) if
$\w(T)$ is an LRS word.

Given three strict partitions $\lambda, \mu, \nu \subset \rho_n$,
define $f(\lambda, \mu ;\nu)$ to be the number of LRS tableaux of
shape $\Sh(\l^\vee/\mu)$ with content $\nu^\vee$.  (If $\mu$ is not
contained in $\l^\vee$ then we set $f(\l, \mu;\nu) = 0$.)  For
example, if $n = 7$ we have $f((5,\!3,\!1),(5,\!2);(6,\!5,\!4,\!1)) =
4$ as counted by the following list of LRS tableaux:
\[ \picB{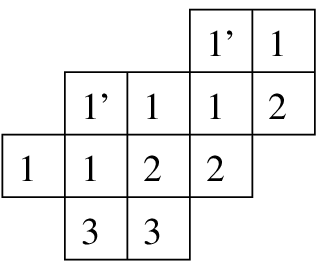} \hspace{6mm} \picB{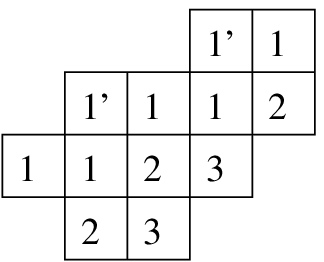} \hspace{6mm}
   \picB{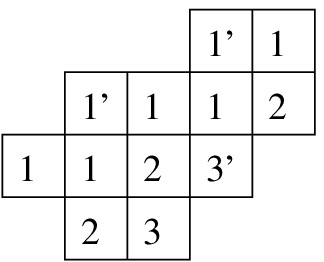} \hspace{6mm} \picB{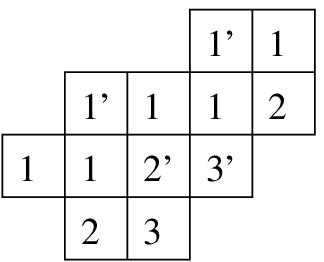} 
\]

\begin{thm}  
\label{typebthm}
The constant $f(\lambda, \mu; \nu)$ is the coefficient of
  $\t_\nu$ in the product $\t_\lambda \cdot \t_\mu$.
\end{thm}
Using the same argument as in the proof of Theorem \ref{typeathm},
Theorem \ref{typebthm} follows from the Pieri rule (\ref{bpieri}) 
and the next proposition, which comes from the associativity relation 
in $H^*(Y,\Z)$.

\begin{prop}
\label{assbprop}  For any three strict partitions 
  $\l,\mu,\nu \subset \rho_n$ and integer
  $p \lequ n$, we have
\begin{equation}
\label{assbeq}
\sum_{\l\xrightarrow{p} 
\widetilde\l} 2^{N(\wt{\l}/\l)}f(\widetilde\lambda, \mu;\nu)
   = \sum_{\mu \xrightarrow{p} 
\widetilde \mu} 2^{N(\wt{\mu}/\mu)}f(\lambda, \widetilde\mu;\nu) \,.
\end{equation}
\end{prop}

The proof of Proposition \ref{assbprop} occupies the remainder of this
section. Define the {\em main diagonal} $\Delta$ to be the set of squares
along the southwest border of $\Sh(\rho_n)$.  We will apply the
shifted analogue of Sch\"utzenberger's sliding operation, constructed
by Worley \cite{W} and Sagan \cite{Sa}, to LRS tableaux.  
This involves the usual sliding moves which refer 
to the alphabet $A$, with the exception of the horizontal slide 
in case (a) below, when a different rule applies. In addition, there
is a {\em special slide} in case (b), which is used only 
when the empty box is on the diagonal $\Delta$.
\[ {\mathrm (a)}\hspace{0.4cm}
   \raisebox{-12pt}{\picB{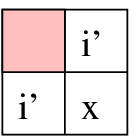}}
   \hspace{5pt}\longmapsto\hspace{5pt}
   \raisebox{-12pt}{\picB{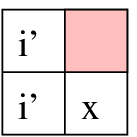}}
   \hspace{1.5cm} {\mathrm (b)}\hspace{0.4cm}
   \raisebox{-12pt}{\picB{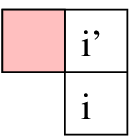}}
   \hspace{5pt}\longmapsto\hspace{5pt}
   \raisebox{-12pt}{\picB{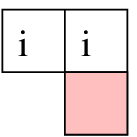}}
\]
These operations are invertible using the obvious reverse slides.

Suppose that we are given an LRS tableau $T$ and an empty box at an
inner corner of $T$, and let $T'$ be the result of performing a
shifted sliding operation to $T$.
The next lemma is parallel to Fact 1, and follows from the fact that
the shifted analogues of the plactic relations preserve the
Littlewood-Richardson-Stembridge property (see \cite{W, Sa, St} for
details.)  We give a direct proof here.

\begin{lemma}
\label{fact1b}
$T$ is an LRS tableau
if and only if $T'$ is an LRS tableau.
\end{lemma}
\begin{proof}
For any $a$ in the alphabet $A$, let
$N_a(\w)$ denote the number of occurrences of $a$ in $\w$. It follows
immediately from the definitions that for any LRS word $\w$, 
\begin{equation}
\label{eqn:decr}
N_i(\w) > N_{i+1}(\w),
 \ \ \ \text{for each unmarked} \  i\in A.
\end{equation}

Since horizontal slides do not change the word of a tableau, we need
only consider special and vertical slides. Observe that in either case, 
condition (ii) in the definition of an LRS tableau is easily verified; 
hence we concentrate on condition (i). 

We start with a special slide as displayed below.

\[ \picB{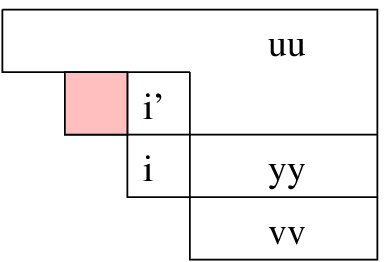}
   \hspace{0.5cm} \raisebox{0.92cm}{$\longmapsto$} \hspace{0.5cm}
   \picB{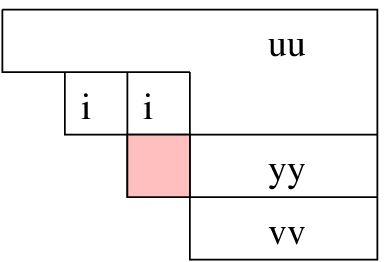}
\]
We must show that if $\w_1=\uu i' \y i \vv$ is an LRS word, then so is
$\w_2=\uu i i \y \vv$. 
Using (\ref{eqn:decr}) we see that $N_i(\uu) +
1 \lequ N_i(\w_1) < N_{i-1}(\w_1) = N_{i-1}(\uu)$.  Since $i',i \not
\in \y \vv \wh{\vv} \wh{\y}$ this implies that every $i'$ and $i$ in
the word $\w_2\wh{\w}_2$ is preceded by more occurrences of $i\m 1$
than of $i$.  Furthermore, since $N_i(\w_1\wh{\vv} (i\p 1)' \wh{\y})
\gequ N_{i+1}(\w_1 \wh{\vv} (i\p 1)' \wh{\y})$ it also follows that
every $(i\p 1)'$ and $i\p 1$ in $\w_2\wh{\w}_2$ is preceded by more
occurrences of $i$ than of $i\p 1$.  All other symbols are not
affected by the slide.

Next, consider a vertical slide.  In the figure, $a$ and $b$ are
symbols such that $a \lequ b$ (if $b$ is marked then $a < b$).
\[ \picB{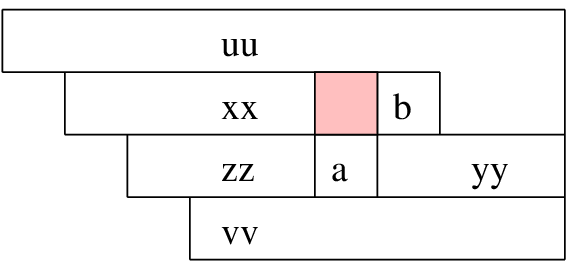} 
   \hspace{0.5cm} \raisebox{0.92cm}{$\longmapsto$} \hspace{0.5cm}
   \picB{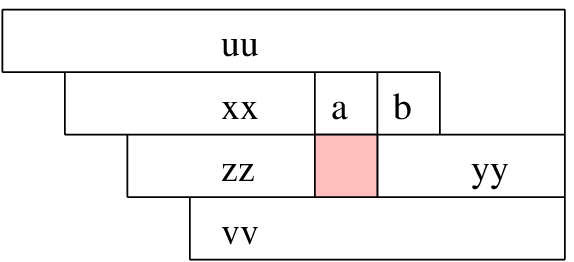}
\]

We must show that if $\w_1=\uu b \x \y a \z \vv$ is an LRS word then
so is $\w_2=\uu b a \x \y \z \vv$.  Assume first that $a=i$ is
unmarked.  To see that every $i'$ and $i$ in the new word $\w_2
\wh{\w}_2$ is preceded by more occurrences of $i\m 1$ than of $i$, we
must show that $N_i(\uu b) + N_{i'}(\x) < N_{i-1}(\uu b)$.  If
$N_i(\z) + N_{i'}(\z) \gequ N_{i-1}(\x) + N_{i'}(\x)$ then this
follows from the LRS condition $N_i(\uu b \x \y i \z) + N_{i'}(\z)
\lequ N_{i-1}(\uu b \x)$.  Otherwise $N_i(\z) + N_{i'}(\z) <
N_{i-1}(\x) + N_{i'}(\x)$ which can only happen when $\z$ is a string
of copies of $i$ terminating at the diagonal $\Delta$, in which case
we have $N_i(\z) = N_{i-1}(\x) + N_{i'}(\x) - 1$ and $i\m 1 \not\in \z
\vv$.  The word $\z$ here cannot contain
$i'$ because of condition (ii) in the definition of an LRS tableau.
Using (\ref{eqn:decr}) we get 
\[
N_i(\uu b) + 1 + N_i(\z) \lequ
N_i(\w_1) < N_{i-1}(\w_1) = N_{i-1}(\uu b) + N_{i-1}(\x)
\]
which also implies the required inequality.

Since $\wh{a} = (i\p 1)'$, we also must check that the string $\w_2
\wh{\vv} \wh{\z} \wh{\y} \wh{\x}$ contains more occurrences of $i$
than of $i\p 1$.  The only way this can fail is if $\wh{\y} \wh{\x}$
contains an $i\p 1$, i.e.\ if $(i\p 1)' \in \x \y$.  Now all symbols in
$\x$ are less than $i$, so $(i\p 1)' \not \in \x$.  If $(i\p 1)' \in
\y$ then $\wh{b} = (i\p 1)'$ or $\wh{b} = i\p 1$, so the lattice
property of the original word $\w_1\wh{\w}_1$ implies the desired one.

Now suppose that $a = i'$ is marked.  To see that the displaced $i'$
is not a problem, we must verify that $N_i(\uu b) < N_{i-1}(\uu)$.
Since $i \not \in \x \y$ and $i\m 1 \not \in \x \y$, 
this follows from the LRS property of the
original word.
We also need to check that all symbols $(i\p 1)'$ and $i\p 1$ in
$\w_2\wh{\w}_2$ are preceded by enough occurrences of $i$.  This can
only fail if $\wh{\y} \wh{\x}$ contains $(i\p 1)'$ or $i\p 1$, i.e.,  
if $\x \y$ contains $i$ or $(i\p 1)'$.  These symbols cannot be in
$\x$ since all symbols in $\x$ are less than $i'$.  The only symbol
among the two that can be in $\y$ is $(i\p 1)'$, and this can only
occur once in $\y$.  Furthermore, we must have $\wh{b} = (i\p 1)'$ or
$\wh{b} = i\p 1$.  Since $i \not\in \wh{\y}\wh{\x}$ and $i\p 1 \in
\wh{\y}$, we deduce that $\w_1 \wh{\vv} \wh{\z}$ contains more
occurrences of $i$ than of $i\p 1$, as required.

By inverting these arguments, one can show that reverse slides also
send LRS tableaux to LRS tableaux. The details are left to the reader.
\end{proof}

\medskip

As in the proof of Proposition \ref{assprop}, we shall call an empty
box contained inside the skew shape $\Sh(\l^\vee/\mu)$ a {\em hole},
but we will need to distinguish between two kinds of holes. For this
purpose, we extend the ordered alphabet $A$ to
$\wt{A}=A\cup\{\hole',\hole\}$, where $\hole'<\hole$ and the new 
symbols represent
a marked and an unmarked hole.  Define a {\em NW-holed tableau}
(respectively, a {\em SE-holed tableau}) to be a filling of a shifted
shape $\Sh(\l^\vee/\mu)$ with symbols from $\wt{A}$ so that the
entries in $A$ satisfy the usual conditions and the holes form a
shifted horizontal strip $L$ along its northwest (respectively,
southeast) border, such that $\w(L)$ is an LRS word. This means that
the holes in a NW-holed tableau occupy a skew shape
$\Sh(\wt{\mu}/\mu)$ for which $\mu\xrightarrow{p} \wt{\mu}$ so that
any hole above another hole is marked, any hole to the right of
another hole is unmarked, and the most southwest hole is unmarked; the
conditions for a SE-holed tableau are similar.

The identity (\ref{assbeq}) is equivalent to the statement that there
are equally many NW-holed and SE-holed LRS tableaux with content $\nu$
on the shape $\Sh(\l^\vee/\mu)$.  We will use shifted slides to
construct an explicit bijection between these two kinds of tableaux.
Given a NW-holed LRS tableau, we first slide the unmarked holes to the
south-east border, in right to left order, after which we slide the
marked holes, proceeding from bottom to top.  If the final position of
an unmarked hole is in a row above the final position of the previous
hole, then we change it to a marked hole.  Marked holes always stay
marked.

\[ \picB{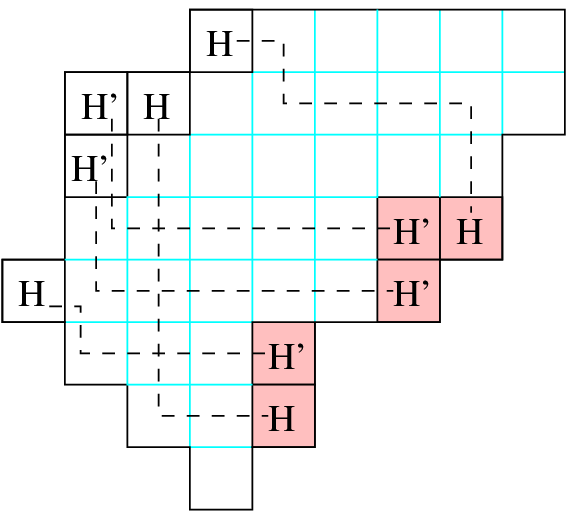} \]

For the reverse bijection, we begin by sliding the marked holes in top
to bottom order, followed by the unmarked holes in left to right
order.  If the path of a marked hole intersects the diagonal $\Delta$
then we erase its marking; the unmarked holes remain unmarked.  To
verify that these two transformations are inverse to each other, we
must check that after all the holes have been slid by one of them,
the other will slide them back in the opposite order.

Let $P$ be a set of boxes in the shifted diagram $\Sh(\l^\vee/\mu)$,
and let $B$ be any box in this diagram. 
We say that $B$ lies {\em west of\/} $P$ if
$P$ contains a box which is strictly east and weakly north of $B$.
And we say that $B$ lies {\em
  north of\/} $P$ if $P$ contains a box which is strictly south and
weakly west of $B$.  

\begin{lemma} \label{lemma:crossb}
Consider the path of a hole $\hole_2$ which slides directly 
after a hole $\hole_1$.
\begin{abcenum}
\item At any given step, if $\hole_2$ lies west of the sliding path 
of $\hole_1$, and $\hole_2$ is not on $\Delta$, then at the next step $\hole_2$
will remain west of the path of $\hole_1$.
\item
At any given step, if $\hole_2$ lies north of the sliding 
path of $\hole_1$,
then the same is true at the next step.
\end{abcenum}
\end{lemma}
\begin{proof}
  Suppose the position of the hole $\hole_2$ is as indicated in the
  figure.
\[ \picB{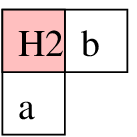} \]

The only way (a) can fail is if $\hole_1$ was in the position of
$b$ and moved down from there.  But then $a \lequ b$ (and if $b$ is marked
then $a < b$), so $\hole_2$ will also move down.  
Notice that there
must be a symbol from $A$ 
in the square occupied by $a$, because $\hole_2$ is
not on the diagonal $\Delta$.

The only way (b) can fail is if the first hole $\hole_1$ was in the position of
$a$ and moved east from there. But this means that $a \gequ b$ (and if $b$
is unmarked then $a > b$), hence $\hole_2$ will move east as well.
This time there must be a symbol from $A$ in the square occupied by $b$.
\end{proof}

\medskip Consider the sequence of slides from northwest to southeast,
beginning with the unmarked holes.
If the path of an unmarked hole crosses the previous path,
then by \reflemma{lemma:crossb} (a) this must be at a corner, and
\reflemma{lemma:crossb} (b) then implies that the hole will remain north
of the previous path from that point onwards.  Since this creates a
path which meets the diagonal $\Delta$, the next unmarked hole will be
forced to cross it, and so on.  The result is that all of the
remaining unmarked holes will become marked and land in 
reverse order. After all the unmarked holes have been slid, 
\reflemma{lemma:crossb} (b) will force every subsequent marked
hole to stay north of the previous hole's path, thus all the marked 
holes retain their order. It follows that
the reverse slides are performed in the opposite order, as
required. Similar arguments can be used to show that reverse slides
will deposit the holes along the northwest border in the opposite
order. This completes the proof of Proposition \ref{assbprop}.

\begin{example}
The following gives an example of the bijection:
\[ \picB{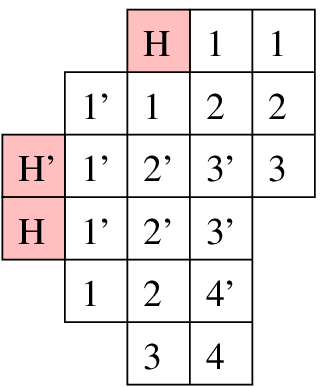} \hspace{2mm} \raisebox{14mm}{$\mapsto$} \hspace{2mm}
   \picB{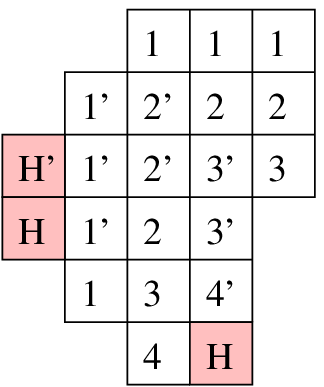} \hspace{2mm} \raisebox{14mm}{$\mapsto$} \hspace{2mm}
   \picB{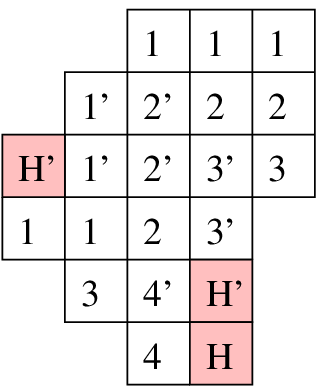} \hspace{2mm} \raisebox{14mm}{$\mapsto$} \hspace{2mm}
   \picB{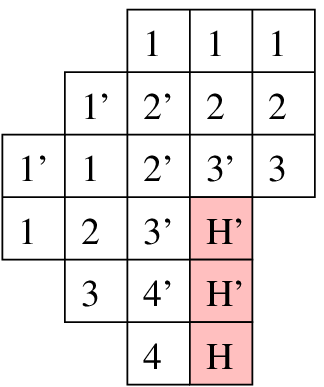} 
\]
\end{example}

\begin{remark}
Arguing as in Section \ref{lrsec}, we can show that
$f(\l,\mu;\nu)$ is equal to the number of LRS tableaux of shape
$\Sh(\nu/\l)$ with content $\mu$, which is Stembridge's original
statement of the rule.  Note also that the even orthogonal
Grassmannian $OG(n+1,2n+2)$ is isomorphic to the odd orthogonal
Grassmannian $OG(n,2n+1)$, and the Schubert structure constants 
for these two spaces coincide. The Schubert classes on the 
Lagrangian Grassmannian $LG(n,2n)$ are also indexed by strict 
partitions $\l$ contained in $\rho_n$, and the corresponding
structure constants $e(\l,\mu;\nu)$ satisfy the identity
$e(\l,\mu;\nu)=2^{\ell(\l)+\ell(\mu)-\ell(\nu)}f(\l,\mu;\nu)$.
This follows by comparing the Pieri formulas for these spaces;
see \cite{Pr} for more details. Therefore, the proof of the
Littlewood-Richardson-Stembridge rule given here also covers 
the maximal isotropic Grassmannians in Lie types $C$ and $D$.
\end{remark}


\end{document}